\documentclass[11pt]{article}
\usepackage{amsmath}
\usepackage{amsfonts}
\usepackage{amssymb}
\usepackage{amsthm}
\usepackage[T1]{fontenc}
\usepackage[latin1]{inputenc}
\usepackage[greek,francais,english]{babel}

\newcommand{\be}{\begin{equation}}
\newcommand{\ee}{\end{equation}}
\newcommand{\ue}{u^{\varepsilon}}
\newcommand{\Se}{S^{\varepsilon}}
\newcommand{\fe}{f^{\varepsilon}}
\newcommand{\me}{m^{\varepsilon}}
\newcommand{\dv}{\mathrm{div}}
\newcommand{\p}{\partial}
\newcommand{\R}{\mathbb R}
\newcommand{\N}{\mathbb N}
\newcommand{\e}{\varepsilon}
\newcommand{\supp}{\mathrm{Supp}}

\newtheorem{thm}{Theorem}[section]
\newtheorem{prop}{Proposition}[section]
\newtheorem{lemma}{Lemma}[section]
\newtheorem{rem}{Remark}[section]

\title{Existence of solutions of the hyperbolic Keller-Segel model}

\author{Beno\^ \i t Perthame\thanks{
D\'epartement de Math\'ematiques et Applications,
\'Ecole Normale Sup\'erieure, CNRS UMR8553 ,
            45 rue d'Ulm, F~75230 Paris cedex 05;
email: perthame@dma.ens.fr}
\and Anne-Laure Dalibard \thanks{
CEREMADE, Universit\'e Paris Dauphine, Place de Lattre de Tassigny, F75775 Paris cedex 16;
email: dalibard@ceremade.dauphine.fr }}

\date{\today}

\begin{document}

\bibliographystyle{amsplain}
\maketitle

\begin{abstract} We are concerned with the hyperbolic Keller-Segel model with quorum sensing,  a model describing the collective  cell movement due to chemical signalling with a  flux limitation for high cell densities.

This is a first order quasilinear equation, its flux
depends on space and time via the solution to an elliptic PDE in
which the right hand side is the solution to the hyperbolic
equation. This model lacks strong compactness or contraction  properties.
Our purpose is to prove the existence of an entropy solution
obtained, as usual, in passing to the limit in a sequence of
solutions to the parabolic approximation.

The method consists in the derivation of a kinetic formulation for
the weak limit. The specific structure of the limiting kinetic
equation allows for a `rigidity theorem' which identifies some
property of the solution (which might be non-unique) to this
kinetic equation. This is enough to deduce a posteriori the strong
convergence of a subsequence.
\end{abstract}


\noindent {\bf AMS Class. No}: 35D05, 35L60,  92C17.
%---------------------------------
\section{Introduction}

We consider the hyperbolic Keller-Segel model \be
\left\{\begin{array}{l} \p_t u + \dv \left(\nabla S(t,y)\; g(u)
\right)=0,\quad t>0, \ y\in \Omega,
\\ \\
u(t=0)=u_0\in L^1\cap L^{\infty}(\Omega),\quad 0\leq u_0\leq 1\ \text{a.e},
\\ \\
-\Delta  S + S =u\quad \text{in}\ \Omega,
\\ \\
\nabla S \cdot n_\Omega=0\quad \text{on} \ \p \Omega.
\end{array}\right.
\label{eq:KShyp} \ee

Here, the function $g(u)$ is given by
$$
g(u)=u (1-u)
$$
therefore we restrict ourselves to solutions satisfying $0 \leq
u(t,x) \leq 1$. The problem is posed on $\Omega$, it is any
bounded domain in $\R^N$, with $\mathcal C^1$ boundary, and
$n_\Omega(y)$ is the outward normal to $\Omega$ at $y\in \p
\Omega$. One can also take the torus $\Omega=\Pi_{i=1}^N (0,T_i)$,
with $T_i>0$, and with periodic boundary conditions in $\Omega$
for $S$; the results and proofs are the same. Notice in particular
that the normal flux in the equation on $u$ vanishes on $\p
\Omega$ and thus the boundary is characteristic; therefore we do
not need boundary conditions for $u$ (this prevents us from
investigating questions which rise specific difficulties, see
\cite{Otto_cras} for instance).

This model represents the density $u(t,y)$ of cells moving with  a
collective chemotactic attraction through the chemical potential
$S$. Their sensitivity is limited by the so-called `quorum
sensing' term $(1-u)$ in $g(u)$. It enters a general class of problems in
the description of cells movement
(\cite{murray,SenbaSuzuki,H03a,H03,HillenP, HillenP2,
BP2,CalvezCarrillo}). Usually a diffusion term is added to
represent the random motion of the cells and the above model
corresponds to the small viscosity limit which has been advocated
by several authors, see \cite{OthmerStevens,
Stevens,DolakSchmeiser,BurgerDolakSchmeiser} for these aspects.

This derivation implies that the system \eqref{eq:KShyp} comes
with an entropy structure as usual (\cite{Serre, Dafermos}). But
extra terms enter in this entropy structure because of the space
dependency of the flux and this leads to  a specific difficulties
(see \cite{BALF, AudusseP, vovelle,kinform} and the references
therein). For any $\mathcal C^2$ convex function $\eta$ (the
so-called entropy), we have \be \frac{\p}{\p t}\eta(u) +
\dv(\nabla S\: q(u)) + (u -S)\left[q-g \eta' \right](u) \leq 0,
\label{in:entrop_hyp} \ee where $q'(\xi):=g'(\xi)\eta'(\xi)$ for
$\xi\in\R$. This accounts for the correct jump condition on
possible discontinuities of $u$. But due to the dependency of the
flux on $\nabla S(t,y)$, the above model poses specific
difficulties compared to the usual theory of quasilinear scalar
conservation laws: no a priori compactness is known in dimension
larger than 1 (no $BV$ bounds or $L^1$ compactness), contraction
principle or uniqueness are not known and averaging lemma for the
kinetic formulation (see below) do not apply (because the
transport is mostly one dimensional in the direction $\nabla S$).
Even time continuity  in $L^1$ does not follow from the method we
develop in this paper. As a consequence we do not know if the full
family of solutions to the diffusion approximation converges, but
only subsequences. All these
questions are left open and seem difficult.\\

Consequently, our proof relies on the weak limit of the diffusion
approximation of \eqref{eq:KShyp} that we study through its
kinetic formulation. Passing to the limit we obtain a weak form of
the kinetic formulation of the hyperbolic limit. The main
ingredient then is to  prove a {\em rigidity theorem} for the
solution which implies that the weak limit is a usual entropy
solution and that subsequences converge strongly. The kinetic
formulation and the main results are presented in the next
subsection. The diffusion limit is studied in section
\ref{sec:parabolic}, and the rigidity theorem is proved in section
\ref{sec:rigidity}. Finally, we  analyze the long time behavior of
solutions in Section \ref{sec:time}. Some technical aspects are
left in an appendix.

%----------------------------------------------------------
\section{Main results}
\label{sec:results}

Our main existence result is the following:

\begin{thm}
The system \eqref{eq:KShyp} has a solution $u \in
L^\infty(\R^+\times \Omega)$, $S\in L^\infty(\R^+;H^1(\Omega))$
satisfying $0\leq u(t,y) \leq 1$, $0\leq S(t,y) \leq 1$, and all
the entropy inequalities  \eqref{in:entrop_hyp}  (in the weak
sense, with initial data $\eta(u_0)$). \label{th:exist}
\end{thm}

Because our method is based on weak limits as mentioned earlier,
it is more convenient to use the kinetic formulation of
\eqref{eq:KShyp} (see \cite{LPTCRAS,LPT,BP,Dalibard} for the
theory of kinetic formulations and recent applications). It is a
way to represent all the inequalities \eqref{in:entrop_hyp} in a
single equation on the unknown defined on
$[0,\infty)\times\Omega\times\R$, $f(t,y,\xi)=\mathbf
1_{\xi<u(t,y)}$, namely \be \left\{\begin{array}{l} \frac{\p f}{\p
t} + (\xi -S)g(\xi)\frac{\p f}{\p \xi} + g'(\xi)\nabla_y S \cdot
\nabla_y f =\frac{\p m}{\p \xi},
\\ \\
m(t,y, \xi)\geq 0 \quad \text{a bounded measure on } [0,T]\times
\Omega \times \R, \quad \forall T>0,
\\ \\
f(0,y, \xi)= \mathbf 1_{\xi<u_0(y)},
\\ \\
-\Delta  S + S =u:=  \int \limits_0^\infty f(t,y,\xi) d\xi \quad
\text{in}\ \Omega,\qquad \quad \nabla S \cdot n_\Omega=0\quad
\text{on} \ \p \Omega.
\end{array} \right.
\label{eq:kinhypint} \ee This is equivalent to
\eqref{in:entrop_hyp}, and one can recover \eqref{in:entrop_hyp}
from \eqref{eq:kinhypint} using that $\eta(u)= \int \eta'(\xi)
\mathbf 1_{\xi<u(t,y)} d\xi$ because we can always take
$\eta(u)=0$ for $u\leq 0$ without loss of generality; see also
Section \ref{sec:parabolic} for an alternative derivation.

The outcome of our proof is the following rigidity theorem:
%-------------------------
\begin{thm} Consider a weak solution to the kinetic equation
\be \left\{\begin{array}{l} \frac{\p f}{\p t} + (\xi
-S)g(\xi)\frac{\p f}{\p \xi} + g'(\xi)\nabla_y S \cdot \nabla_y f
 + R(t,y,\xi)=\frac{\p m}{\p \xi},
\\ \\
m(t,y, \xi)\geq 0 \quad \text{a bounded measure on } [0,T]\times
\Omega \times \R, \quad \forall T>0,
\\ \\
f(0,y, \xi)= \mathbf 1_{\xi<u_0(y)},
\\ \\
-\Delta  S + S =u:=  \int \limits_0^\infty f(t,y,\xi) d\xi \quad
\text{in}\ \Omega,\qquad \quad \nabla S \cdot n_\Omega=0\quad
\text{on} \ \p \Omega.
\end{array} \right.
\label{eq:kinhyprem} \ee which satisfies the properties
\\
(i) $\; 0\leq f(t,y, \xi)\leq 1$ and $f=1$ for $\xi < 0$, $f=0$
for $\xi >1$, $f$ is nonincreasing in $\xi$,
\\
(ii) there exists a constant $C>0$ such that $|R|\leq C f(1-f)$
almost everywhere,
\\
(iii) the measure $\; m$ vanishes for $\xi <0$ or $\xi >1$.
\\
Then, we have $f(t,y,\xi)=\mathbf 1_{\xi<u(t,y)}$ and $u(t,y)$ is
an entropy solution to  \eqref{eq:KShyp}. \label{th:kin}
\end{thm}

The proof of these two results is given in the next sections. The
strategy in the following : as in \cite{DolakSchmeiser}, we take a
parabolic approximation of \eqref{eq:KShyp}, and we intend to pass
to the limit in its solution $u^\e$ as the viscosity vanishes.
However, unlike in \cite{DolakSchmeiser} and as mentioned earlier,
the problem \eqref{eq:KShyp} lacks a priori compactness bounds for
$u^\e $ when $N>1$. Hence, we rather pass to the limit in a
kinetic formulation of the approximate problem. The weak limit of
the sequence $\fe=\mathbf 1_{\xi<\ue(t,y)}$ satisfies equation
\eqref{eq:kinhyprem}, with a remainder $R$ which can be explicitly
computed in terms of $f$ and which satisfies (ii). Thus theorem
\ref{th:kin} implies in turn that $\ue$ converges strongly to $u$.
Let us finally mention that an alternative proof for the local
existence of strong solutions can be carried out, see
\cite{BurgerDolakSchmeiser}; however, as stressed by M. Burger, Y.
Dolak and C. Schmeiser in \cite{BurgerDolakSchmeiser}, their
strategy does not yield any information on the global existence of
weak solutions.

%----------------------------------------------------------
\section{The parabolic limit}
\label{sec:parabolic}

In this section, we introduce and study the approximate parabolic system with $\e>0$:
 \be \left\{
\begin{array}{l}\p_t{\ue}+ \dv \left(\nabla \Se(t,y)\:
\ue (1-\ue) \right)-\e \Delta \ue=0,\quad t>0, \ y\in \Omega\\
\ue(t=0)=u_0\in L^1\cap L^{\infty}(\Omega),\quad 0\leq u_0\leq 1\ \text{a.e},\\
-\Delta \Se + \Se =\ue\quad \text{in}\ \Omega, \\
\nabla \Se \cdot n_\Omega=0\quad \text{on} \ \p \Omega,\\
\nabla \ue \cdot n_\Omega=0\quad \text{on} \ \p \Omega\quad
\text{for a.e.}\ t>0.
\end{array}\right.
\label{eq:KSpara} \ee

Our goal is to pass to the weak limit in this system but we first state the following result
\begin{prop}
There exists a unique solution $(\ue,\Se)\in
L^2_{\text{loc}}(0,\infty; H^1(\Omega))\times L^{\infty}(0,\infty;
H^1(\Omega))$ of the problem \eqref{eq:KSpara} and it  satisfies
the following bounds : for all $1\leq q<\infty$, for all $T>0$,
there exist constants $C_1(N, \Omega,q)$, $C_2(N,\Omega,T)$ such
that
\begin{gather}
0\leq \ue(t,y)\leq 1\quad \text{a.e. on}\ [0,\infty)\times \Omega,\label{bound1}\\
0\leq \Se(t,y)\leq 1\quad \text{a.e. on}\ [0,\infty)\times \Omega,\label{bound1b}\\
||\Se||_{L^{\infty}(0,\infty; W^{2,q}(\Omega))}\leq C_1,\label{bound2}\\
\sqrt{\e} ||\nabla \ue||_{L^2((0,T)\times \Omega)} + ||\p_t
\Se||_{L^2((0,T); H^1(\Omega))}\leq C_2.\label{bound3}
\end{gather}
And for any $\mathcal C^2$ convex function $\eta$, we have with the notation in  \eqref{in:entrop_hyp},
\be
\frac{\p}{\p t}\eta(\ue) + \dv(\nabla \Se\: q(\ue)) + (\ue -
\Se)\left[q-g \eta' \right](\ue) -\e \Delta \eta(\ue)\leq 0.
\label{in:entrop_para} \ee
\label{prop:epsmodel}
\end{prop}

\begin{proof}
Existence and uniqueness of $(\ue,\Se)$ are easily proved thanks
to semi-group techniques. The bounds \eqref{bound1} follows from
the maximum principle because $0$ and $1$ are solutions for all
drifts $\nabla S^\e$, the bound \eqref{bound1b} also follows from
the maximum principle, whereas \eqref{bound2} is the regularizing
effect for elliptic equations with smooth coefficients.

The first bound on $\nabla \ue$ in \eqref{bound3} is obtained by
multiplying by $\ue$ the evolution equation on $\ue$. Eventually,
differentiating the equation giving $\Se$ with respect to $t$
gives
$$
-\Delta \p_t \Se + \p_t \Se = - \dv \left( \nabla\Se
\:\ue(1-\ue)\right) + \e \Delta \ue,
$$
and the right-hand side is bounded in $L^2_{\text{loc}}(0,\infty;
H^{-1}(\Omega))$ uniformly in $\e$; the second bound of
\eqref{bound3} follows.

The entropy inequality \eqref{in:entrop_para} is obtained by
multiplication of the evolution equation by $\eta'(\ue)$ and using
the chain rule.
\end{proof}

Next, we pass to the limit in the system
\eqref{eq:KSpara}. However, because they do not provide strong compactness, the bounds on the
sequence $\ue$ are insufficient to pass to the limit in the
nonlinear term
$$
\nabla \Se \ue (1-\ue).
$$
In \cite{DolakSchmeiser}, for $N=1$, strong compactness is
obtained thanks to uniform $BV$ bounds on the sequence $\ue$;
however, as we have already pointed out, such bounds no longer
hold when $N>1$. Consequently, we pass to the (weak) limit in the
kinetic formulation for problem \eqref{eq:KSpara}. Our next goal
is to present this argument.

We take $\eta(u)=(u-\xi)_+$ in \eqref{in:entrop_para}, with $\xi
\in \R$, and we differentiate (in the distributional sense) the
inequality obtained with respect to $\xi$. This yields
 \be \frac{\p \fe}{\p t} + (\xi -
\Se)g(\xi)\frac{\p \fe}{\p \xi} + g'(\xi)\nabla_y \Se\cdot
\nabla_y \fe  -\e \Delta_y \fe=\frac{\p \me}{\p \xi},
\label{eq:kinpara} \ee
where
$\me(t,y,\xi)$ is a nonnegative measure on $[0,\infty)\times
\Omega \times \R$. It can be written  explicitely in terms of $\ue$, namely
\begin{eqnarray*}
\me(t,y,\xi)&:=&-\left\{\p_t(\ue-\xi)_+ +
\dv\left(\nabla \Se\:\mathbf 1_{\xi<\ue}(g(\ue)-g(\xi))\right)\right.\\
&& \qquad- \left.(\ue - \Se)\mathbf 1_{\xi<\ue}g(\ue) -\e \Delta_y
(\ue-\xi)_+ \right\}
\\&=&\e |\nabla_y\ue(t,y)|^2 \delta(\xi=\ue(t,y)).
\end{eqnarray*}

Notice that $0\leq \fe\leq 1$ almost everywhere, and
$\fe(t,y,\xi)=0$ when $\xi >1$, $\fe(t,y,\xi)=1$ when $\xi<0$.
Moreover, $\me(t,y,\xi)=0$ when $\xi<0$ or $\xi>1$ (in the sense
of distributions), and $\{\me(t,y,\xi)\}_{\e>0}$ is a family of
 bounded measures on $[0,T]\times \Omega \times \R$, $\forall T>0$.

Hence, there exists a subsequence, still denoted by $\e$, and
functions $u=u(t,y)\in L^{\infty}((0,\infty)\times \Omega)$,
$f=f(t,y,\xi)\in L^{\infty}((0,\infty)\times \Omega\times\R)$,
$S=S(t,y)\in L^{\infty}(0,\infty; W^{2,q}(\Omega))$, and a
nonnegative measure $m=m(t,y,\xi)$ such that, locally in time,
\begin{gather*}
\ue \rightharpoonup u \quad w^* -L^{\infty},\\
\fe \rightharpoonup f \quad w^* -L^{\infty},\\
\me \rightharpoonup m \quad w-M^1,\\
\Se \rightarrow S \quad \text{in}\ L^p_{\text{loc}}(0,\infty;
W^{1,p}(\Omega))
\end{gather*}
for all $p$, $1\leq p < \infty$.

Thus, we can pass to the limit as $\e\to 0$ in equation
\eqref{eq:KSpara}.  All the terms can pass to the limit because they are written as  weak-strong products except $
g'(\xi)\nabla_y \Se\cdot \nabla_y \fe$ which yields an extra term. Indeed, we can write
\begin{eqnarray*}
g'(\xi)\nabla_y \Se\cdot \nabla_y \fe&=&\dv_y\left(
g'(\xi)\nabla_y
\Se\;\fe\right) - g'(\xi)\Delta_y \Se\; \fe\\
&=&\dv_y\left( g'(\xi)\nabla_y \Se\;\fe\right) +
(\ue-\Se)g'(\xi)\fe.
\end{eqnarray*}
In the sense of distributions, as $\e\to 0$, we have
\begin{gather*}
\dv_y\left( g'(\xi)\nabla_y \Se\;\fe\right)\rightharpoonup
\dv_y\left( g'(\xi)\nabla_y S\;f\right),\\
\Se\;g'(\xi)\fe\rightharpoonup S\: g'(\xi)\:f.
\end{gather*}
But at this stage, we cannot assert that the weak limit of
$\ue\:\fe$ is $uf$ (but it is possible to identify it, see
\eqref{eq:rolimit} below). Nevertheless, we know that
$\{\ue\fe\}_{\e>0}$ is bounded in $L^{\infty}$; thus, extracting a
further subsequence if necessary, there exists a function
$\rho=\rho(t,y,\xi)$ such that \be \ue\:\fe\rightharpoonup
\rho\quad w^* -\ L^{\infty}. \ee Consequently,
$$
g'(\xi)\nabla_y \Se\cdot \nabla_y \fe\rightharpoonup
g'(\xi)\nabla_y S\cdot \nabla_y f +g'(\xi)\left(\rho-u\:f \right),
$$
and $f$ is a solution of \begin{gather} \p_t f + (\xi - S)g(\xi)
\p_{\xi} f + g'(\xi) \nabla_y S \cdot \nabla_y f +
g'(\xi)(\rho-uf)=\p_{\xi }m,\nonumber\\
-\Delta S + S=u(t,y) \quad \text{in}\ \Omega, \label{eq:kinhyp}\\
\nabla_y S \cdot \nabla n_\Omega=0\quad \text{on}\ \p \Omega,\nonumber\\
f(t=0,y,\xi)=\mathbf 1_{\xi<u_0(y)}.\nonumber
\end{gather}

Moreover, $f$, $u$ and $m$ inherit the following properties
\begin{gather*}
0\leq f\leq 1\quad \text{a.e.},\\ f(t,y,\xi)=0\quad \text{when}\
\xi>1,\quad  f(t,y,\xi)=1\quad \text{when}\ \xi<0,\\
m(t,x,\xi)=0\quad \text{when}\ \xi>1\ \text{or }\xi<0,\\
\int_0^T \int_{\Omega}\int_{\R} m(t,y,\xi)\;
dt\:dy\:d\xi<\infty\quad \forall T>0.
\end{gather*}
And there exists a nonnegative measure $\nu(t,y,\xi)$ such that
$\nu((0,T)\times \bar{\Omega}\times \R)<\infty$ for all $T>0$ and
\be \p_{\xi} f(t,x,\xi)=-\nu(t,x,\xi)\leq 0\label{def:nu} \ee in
the sense of distributions. This follows from the fact that
$$
\p_{\xi }\fe(t,y,\xi)=-\delta(\xi-\ue(t,y)).
$$
This means that we have derived the properties (i) and (iii)
assumed in Theorem \ref{th:kin}. The remainder term $R$ is here
equal to $R(t,y,\xi):=g'(\xi)(\rho-u f)(t,y,\xi)$. There remains
to derive a formula for $\rho$ which we do now.

Let
$\varphi_1\in \mathcal C^{\infty}_0([0,\infty)\times \Omega)$,
$\varphi_2\in \mathcal C^{\infty}_0(\R)$ be test functions.

Then
\begin{eqnarray*}
&&\int_0^\infty\int_{\Omega\times\R}\rho(t,y,\xi)\varphi_1(t,y)\:\varphi_2'(\xi)\:dt\:dy\:d\xi
\\&=& \lim_{\e\to
0}\int_0^\infty\int_{\Omega\times\R}\ue(t,x)\fe(t,x,\xi)\varphi_2'(\xi)\varphi_1(t,y)\:dt\:dy\:d\xi\\
&=&\lim_{\e\to
0}\int_0^\infty\int_{\Omega}\ue(t,x)(\varphi_2(\ue(t,x)))\varphi_1(t,y)\:dt\:dy\\
&=&\lim_{\e\to
0}\int_0^\infty\int_{\Omega\times\R}\frac{d}{d\xi}(\xi
\varphi_2(\xi))\fe(t,x,\xi)\varphi_1(t,y)\:dt\:dy\:d\xi\\
&=&\int_0^\infty\int_{\Omega\times\R}\frac{d}{d\xi}(\xi
\varphi_2(\xi))f(t,x,\xi)\varphi_1(t,y)\:dt\:dy\:d\xi\\
&=&\int_0^\infty\int_{\Omega\times\R}\frac{d}{d\xi}(\xi
\varphi_2(\xi))f(t,x,\xi)\varphi_1(t,y)\:dt\:dy\:d\xi
\end{eqnarray*}
Consequently,
$$
-\frac{\p}{\p \xi}\left[\rho -\xi f \right]= f.
$$
Next, we integrate this equation on $\R$ ($t,y$ are treated as
fixed parameters), with the boundary conditions $f(t,y,\xi)=0$ and
$\rho(t,y,\xi)=0$ when $\xi>1$. We get
\be
\rho(t,y,\xi) -\xi f(t,y,\xi) = \int_{\xi}^{\infty}
f(t,y,\xi')\:d\xi'.
\label{eq:rolimit}
\ee

We show later (see lemma \ref{lem:rho_uf}) that this implies the
assumption (ii) on $R$ in theorem \ref{th:kin}.

At this stage we have derived the full kinetic formulation for our
problem, which means that the assumptions of Theorem \ref{th:kin}
have been obtained in the (weak) limit of solutions to the
parabolic equation \eqref{eq:KSpara}. We can turn to its proof.

%-------------------------------------------------------
\section{Proof of the rigidity Theorem  \ref{th:kin}}
\label{sec:rigidity}

The technique introduced in \cite{artBP}  is then to compare $f$
and $f^2$ in order to prove that $f$ only takes the values 0 and 1
almost everywhere. Thanks to the monotony assumption in (i) (see
\eqref{def:nu}), we can then deduce easily that there exists
$u=u(t,y)$ such that $f(t,x,\xi)=\mathbf 1_{\xi<u(t,y)}$.
\\

Hence, we multiply \eqref{eq:kinhyprem} by $2f$ and we formally
derive an equation for $f^2$; the difference $f-f^2$ satisfies
\begin{multline} \frac{\p}{\p t}(f-f^2)  + (\xi - S)g(\xi) \p_{\xi}(f-f^2) +
g'(\xi) \nabla_y S \cdot \nabla_y (f-f^2)  + R(1-2f)=\\=\p_{\xi
}m(1-2f). \label{eq:difference}\end{multline} We emphasize that
this calculation, and the following, seems entirely formal;
indeed, since $f$ is not smooth,  the chain rule $2\p_t f f =\p_t
f^2$ for instance, has to be justified. Thus, regularizations in
$(t,y,\xi)$ are necessary in order to make the argument rigorous.
Those are fairly standard (see \cite{DPL,BP,artBP}), and will be
detailed in the Appendix.

It can be seen in the above equation that the key of our method is
the assumption (ii) on the term $R$. In the case when $R$ is equal
to $R=g'(\rho-uf)$, with $\rho$ given by \eqref{eq:rolimit}, the
inequality in assumption (ii) is proved in the following
%----------------------------
\begin{lemma} For $T>0$, set
$$
C:=\limsup ||\ue||_{L^{\infty}}((0,T)\times \Omega).
$$
(Notice that $C\leq 1$ here). Then, with $\rho$ given in
\eqref{eq:rolimit}, we have
$$
\left|\rho(t,y,\xi)-u(t,y)\:f(t,y,\xi)\right| \leq C
f(t,y,\xi)\left(1 - f(t,y,\xi) \right)
$$
for a.e. $t\in(0,T)$, $y\in \Omega$, $\xi\in\R$.\label{lem:rho_uf}
\end{lemma}
%--------------------------------
\begin{proof}
From \eqref{eq:rolimit}, for almost every
$(t,y,\xi)\in[0,\infty)\times \Omega \times\R$,

\begin{eqnarray*}
&&\rho(t,y,\xi)-u(t,y)f(t,y,\xi)\\&=&\xi f(t,y,\xi)
+\int_{\xi}^{\infty} f(t,y,\xi')\:d\xi'-u(t,y)f(t,y,\xi)\\
&=&\int_0^{\xi} d\xi'\; f(t,y,\xi) +\int_{\xi}^{\infty}
f(t,y,\xi')\:d\xi'\\&&-f(t,y,\xi)\int_0^{\xi}f(t,y,\xi')\:d\xi' -
f(t,y,\xi)\int_{\xi}^{\infty}f(t,y,\xi')\:d\xi'\\
&=&\left[\int_0^{\xi}(1-f(t,y,\xi')) d\xi'\right] f(t,y,\xi)
+\left[\int_{\xi}^{\infty}
f(t,y,\xi')\:d\xi'\right](1-f(t,y,\xi)).
\end{eqnarray*}
Now, remember that $f(t,y,\xi)$ is decreasing with respect to
$\xi$ (recall \eqref{def:nu}). Therefore $f(t,y,\xi')\leq
f(t,y,\xi)$ for $\xi'\geq \xi$, and $1-f(t,y,\xi')\leq
1-f(t,y,\xi)$ for $\xi'\leq \xi$. And for $t\in[0,T]$, $y\in
\Omega$, $f(t,y,\xi)=0$ for $\xi>\limsup
||\ue||_{L^{\infty}}((0,T)\times \Omega)$.

Eventually, we obtain
$$
0\leq \rho(t,y,\xi)-u(t,y)f(t,y,\xi)\leq\limsup
||\ue||_{L^{\infty}} \left[f(1-f)\right](t,y,\xi).
$$

\end{proof}

Now, we integrate \eqref{eq:difference} on $\Omega\times\R$
(notice that for $\xi<0$ or $\xi>1$, $f=f^2$). We get
\begin{eqnarray}
\frac{d}{dt}\int_{\Omega\times \R}(f-f^2)&\leq &
\int_{\Omega\times\R}(f-f^2)\left\{\frac{\p}{\p\xi}\left[(\xi-S)
g(\xi)
\right] + \Delta_y S g'(\xi) \right\}\nonumber\\
&& + C\int_{\Omega\times\R}\left| (1-2f) \right|(f-f^2)\nonumber\\
&& + 2 \int_{\Omega\times\R}m(t,y,\xi)\:\p_{\xi}f(t,y,\xi)\;dy\:d\xi\nonumber\\
&\leq & C\int_{\Omega\times\R}(f-f^2).\label{in:comp_f/f2}
\end{eqnarray}
In the above inequality, we have used the facts that $f-f^2\geq 0$
because $0\leq f\leq 1$, and $0\leq S\leq 1$ by the maximum
principle. Also,
\begin{gather*}
\left|\Delta S\right|=|S-u|\leq 1,\\
|g(\xi)|,|g'(\xi)|\leq 1\quad \forall \xi\in[0,1],\\
\int_{\Omega\times\R}m(t,y,\xi)\:\p_{\xi}f(t,y,\xi)\;dy\:d\xi=-\int_{\Omega\times\R}m(t,y,\xi)\nu(t,y,\xi)\;dy\:d\xi\leq
0.
\end{gather*}

Consequently, by Gronwall's lemma, we get
$$
0\leq \int_{\Omega\times\R}(f-f^2)(t,y,\xi)\;dy\:d\xi\leq e^{C
t}\int_{\Omega\times\R}(f-f^2)(t=0,y,\xi)\;dy\:d\xi
$$
But $f(t=0,y,\xi)=\mathbf 1_{\xi<u_0(y)}$, and thus
$(f-f^2)(t=0)=0$. We deduce that $f(t,y,\xi)=f^2(t,y,\xi)$ for
a.e. $(t,y,\xi)$, and $f=0$ or $f=1$ almost everywhere. Since $f$
is decreasing in $\xi$, $f=\mathbf 1_{\xi<u(t,y)}$, and it is
easily checked that in that case, $\rho=uf$. Hence $f$ is a
solution of
$$
\p_t f + (\xi-S) g(\xi) \p_{\xi} f + \nabla_y S \cdot \nabla_y f
g'(\xi)=\p_{\xi} m,
$$
and $u$ is an entropy solution of
$$
\p_t u + \dv_y (\nabla S g(u))=0,\quad t>0,\ y\in \Omega.
$$

\begin{rem}
It can be checked that
$$
\mathbf 1_{\xi<\ue(t,y)}\rightharpoonup\mathbf 1_{\xi<u(t,y)}\
w^*- L^{\infty}\quad \iff \quad \ue \to u \quad \text{in}\
L^1_{\text{loc}}((0,\infty)\times \Omega).
$$

Hence, we have proved here the following more general result : let
$\{(u_n, S_n)\}_{n\geq 0}$ be a sequence such that $f_n:=\mathbf
1_{\xi<u_n(t,y)}$ satisfies
\begin{gather*}
\p_t f_n + (\xi-S_n)g(\xi)\p_{\xi} f_n +g'(\xi)\nabla S_n
\cdot\nabla f_n =\p_{\xi} m_n + r_n,\\
-\Delta S_n + S_n =u_n,\\
0\leq u_n\leq 1,\\
f_n(t=0)=1_{\xi<u^0_n(y)},
\end{gather*}
with $m_n$ a nonnegative measure and $r_n\rightharpoonup 0$ in the
sense of distributions.

Assume that $f_n\rightharpoonup f(t,y,\xi)$ $w^*- L^{\infty}$ and
$u^0_n(y)\rightarrow u^0(y)$ strongly in $L^1(\Omega)$ as
$n\to\infty$.

Then there exists $u=u(t,y)\in L^{\infty}((0,\infty)\times
\Omega)$ such that $f=\mathbf 1_{\xi<u(t,y)}$ and $u_n\to u$ in
$L^1((0,T)\times \Omega)$ for all $T>0$. And $u$ is an entropy
solution of \eqref{eq:KShyp}.
\end{rem}

%--------------------------------------------
\section{Long-time behavior}
\label{sec:time}

We wish to mention here a few simple facts on the long-time
behavior of a solution $(u,S)$ of system \eqref{eq:KShyp}. Our
motivation comes from the unusual complexity of the behavior
exhibited in \cite{BurgerDolakSchmeiser,DolakSchmeiser} for this
limit. From their study there appears to be differences between
the large time dynamics of the parabolic system \eqref{eq:KSpara}
and the hyperbolic system \eqref{eq:KShyp}. The
 stability of the steady states of
\eqref{eq:KShyp} is also discussed in these references. The
numerical simulations and formal computations presented in these
two papers also convey a rather good insight of the long time
behavior of solutions. Precisely, the numerical simulations
indicate that for the hyperbolic system \eqref{eq:KShyp} in
dimension one, solutions converge to piecewise constant steady
states as time goes to infinity. In these steady states, regions
of vacuum ($u_\infty=0$ below) are separated from regions where cells
aggregate ($u_\infty=1$ below) by entropic shocks. On the contrary, in
the parabolic system \eqref{eq:KSpara}, when the parameter $\e$ is
large enough ($\e>\frac 1 4$), solutions converge to the constant
solution. And when the diffusivity parameter $\e$ is small, a
metastable behavior occurs : solutions first get close to the
piecewise constant steady states of the hyperbolic system, and
then the regions of cell aggregates (called `plateaus') move
slowly and at last merge with one another.

This section aims at proving that any entropy solution to
\eqref{eq:KShyp}, as built in Theorem \ref{th:exist}, converges
(in a sense detailed in the next proposition) to a steady state
solution: \be \left\{
\begin{array}{l} \dv \left( \nabla S_\infty \; u_\infty
(1-u_\infty)\right)=0\quad\text{in}\ \Omega,
\\ \\
-\Delta S_\infty + S_\infty=u_\infty\quad\text{in}\ \Omega,
\\ \\
\nabla S_\infty \cdot n_{\Omega}=0\quad \text{on}\ \p\Omega,
\end{array}\right.
\label{eq:hypkssta} \ee

Our analysis relies on the energy dissipation inherited from the natural free energy structure for the chemotaxis systems (\cite{CalvezCarrillo,BP2,H03})
\be
\begin{array}{rl}
\frac{d}{dt}\int_{\Omega} u(t) S(t) &=\frac{d}{dt} ||S(t)||_{H^1(\Omega)}^2
\\ \\
&=2\int_{\Omega}|\nabla S(t,y)|^2\:u(t,y)(1-u(t,y))\:dy\geq 0,
\end{array}
\label{eq:endiss}\ee
and consequently.
\be \int_{\tau}^{\infty}\int_{\Omega}|\nabla
S(t,y)|^2\:u(t,y)(1-u(t,y))\:dy\:dt\to 0\quad \text{as}\ \tau
\to\infty.\label{cvgto0}
\ee

The equality \eqref{eq:endiss} is proved in \cite{DolakSchmeiser}
when $N=1$, but the dimension does not play a significant part
here; we reproduce a short proof for the reader convenience.
First, notice that
$$
-\Delta S + S =u, \qquad -\Delta \p_tS + \p_tS =\p_tu,
$$
and thus, multiplying the first equation by $\p_t S$ and the second by $S$, after integration by parts we obtain
\begin{eqnarray*}
\int_{\Omega} u \p_t S &=& \int_{\Omega} \left(\nabla S\cdot
\nabla\p_t S + S\:\p_t S\right)\\
&=&\int_{\Omega}\p_t u S =\frac 1
2\frac{d}{dt}||S||_{H^1(\Omega)}^2.
\end{eqnarray*}
Notice also that $\int_\Omega u S=||S||_{H^1}^2$. Consequently,
\begin{eqnarray*}
\frac{d}{dt}\int_{\Omega} u S &=&2\int_{\Omega} S \p_t u  \\
&=&-2\int_{\Omega}S \; \dv\left(\nabla S g(u) \right)
=2\int_{\Omega}\left|\nabla S\right|^2 g(u).
\end{eqnarray*}
Therefore, we are led to
$$
\int_{\Omega} u(t,y)\frac{\p S(t,y)}{\p t}\:dy =\int_{\Omega}
S(t,y)\frac{\p u(t,y)}{\p t}\:dy=\frac 1 2
\frac{d}{dt}||S||_{H^1(\Omega)}^2.
$$
And \eqref{eq:endiss} follows.

Integrating this equality from $t=0$ to $t=T$, we obtain
$$
\int_0^T \int_{\Omega}|\nabla S|^2 g(u) = \frac 1
2\left[\int_{\Omega}\left(u(T) S(T) - u(t=0)
S(t=0)\right)\right]\leq \frac 1 2 |\Omega|.
$$
Thus
$$
\int_0^\infty\int_\Omega|\nabla S|^2 g(u)\leq \frac 1 2
|\Omega|<+\infty
$$
and \eqref{cvgto0} follows.
\\

We can now state our main result
%---------------------------------
\begin{prop}
Let $(u,S)$ be a global weak solution of
\eqref{eq:KShyp} as in Theorem \ref{th:exist}. Then, for $k\in \N$, let $u_k(t,y):=u(t+k,y)$,
$S_k(t,y):=S(t+k,y)$, $t>0, \ y\in\Omega$.
Then there exists a subsequence $(n_k)_{k\in \N}$
such as when  $k\to\infty$,
\begin{gather}
u_{n_k}(t,y)\rightharpoonup u_{\infty}(y)\quad
w^*-L^{\infty}(\R^+\times\Omega), \label{time1}
\\
S_{n_k}(t,y)\rightarrow S_{\infty}(y)\quad \text{in} \ L^2(0,T;
H^1(\Omega)), \qquad \forall T>0, \label{time2} \\
\int_0^T\int_{\Omega}\mathbf 1_{\nabla S \neq 0}|u_{n_k}(t,y)
-u_{\infty}(y)|\:dt\:dy\to 0,\label{strongcvg}
\end{gather}
where $u_{\infty}=u_{\infty}(y)\in L^{\infty}(\Omega)$,
$S_{\infty}=S_{\infty}(y)\in H^2(\Omega)$ are solutions to
\eqref{eq:hypkssta} and $0\leq  u_{\infty}(y) \leq1$, $0\leq
S_{\infty}(y) \leq1$ and $|\nabla S_\infty| u_\infty(1-u_\infty)=0$. \label{prop:lgtime}
\end{prop}

\begin{proof}

\noindent{\it First step. Weak convergence of $(u_k,S_k)$.} The bounds $0\leq u_k,S_k \leq 1$, and
the elliptic regularity
$$
||S_{n_k}||_{W^{1,2}((0,T)\times \Omega)} + ||\nabla
S_{n_k}||_{W^{1,2}((0,T)\times \Omega)^N}\leq C\quad \forall
k\in\N
$$
provide us directly with \eqref{time1},  \eqref{time2} after extracting subsequences.

\vskip1mm

\noindent{\it Second step. The limits $u_{\infty}$ and
$S_{\infty}$ are independent of time.}

First, since the couple $(u_{n_k}, S_{n_k})$ is a solution of
\eqref{eq:KShyp}, considering a test function $\varphi=\varphi(t,y)$  in
$\mathcal C^{\infty}_0((0,T)\times Y)$, with
$\varphi(0,y)=\varphi(T,y)=0$ for all $y\in\Omega$, then
$$
\int_0^T\int_\Omega u_{n_{k}}\p_t \varphi = -
\int_0^T\int_\Omega\nabla S_{n_k}\cdot\nabla \varphi g(u_{n_k}).
$$
The right-hand side goes to 0 as $k\to\infty$ according to
\eqref{cvgto0} for any test function $\varphi$. Thus
$$
\int_0^T\int_\Omega u_{\infty}(t,y)\p_t \varphi(t,y)\:dt\:dy=0
$$
for any test function $\varphi$ vanishing at $t=0$ and $t=T$.
Consequently, $u_{\infty}$ is independent of $t$ :
$u_{\infty}(t,y)=u_{\infty}(y)$.

Furthermore, in the weak limit it holds
$$
-\Delta S_{\infty} + S_{\infty}=u_{\infty}\quad \text{on}
(0,T)\times\Omega, \qquad \nabla S_{\infty} \cdot n_\Omega=0 \quad
\text{on} \ \p \Omega.
$$
and by uniqueness for this problem, we also have
$S_{\infty}=S_{\infty}(y)$.

\vskip1mm

\noindent{\it Third step. The limiting equation on $u_\infty$.}

We now introduce the notations
\begin{gather*}
A:=\{(t,y)\in(0,T)\times \Omega ;\ \nabla S_{n_k}(t,y)\to 0 \},\\
B:=\{(t,y)\in(0,T)\times \Omega ;\  u_{n_k}(t,y)\to 0 \},\\
C:=\{(t,y)\in(0,T)\times \Omega ;\  u_{n_k}(t,y)\to 1 \}.
\end{gather*}
Then $\lambda\left((0,T)\times \Omega\setminus (A\cup B\cup C)
\right)=0$, where $\lambda$ is the Lebesgue measure, because from \eqref{cvgto0}, we deduce that
$$
\int_0^T \int_{\Omega}|\nabla S_k|^2(t,y) \:u_k(t,y)
(1-u_k(t,y))\:dt\:dy\to 0
$$
as $k\to\infty$. Hence there exists a subsequence, still denoted by
$n_k$, such that
$$
|\nabla S_{n_k}|^2(t,y) \:u_{n_k}(t,y)(1-u_{n_k}(t,y))\to 0\quad \text{a.e.}
$$

It follows
from the above strong convergence results that
\begin{gather*}
\int_0^T \int_\Omega\mathbf 1_A \left|\nabla S_{n_k} \right|^2\to
0=\int_0^T \int_\Omega\mathbf 1_A \left|\nabla
S_{\infty}\right|^2,\\
\int_0^T \int_\Omega\mathbf 1_B u_{n_k} \to 0=\int_0^T
\int_\Omega\mathbf 1_B u_{\infty},\\
\int_0^T \int_\Omega\mathbf 1_C (1-u_{n_k})\to 0=\int_0^T
\int_\Omega\mathbf 1_C (1-u_{\infty}).
\end{gather*}
Consequently,
$$
\int_0^T\int_{\Omega}|\nabla S_{\infty}|^2g(u_{\infty})= T
\int_{\Omega}|\nabla S_{\infty}|^2g(u_{\infty})=0,
$$
and $|\nabla S_{\infty}|^2g(u_{\infty})=0$ almost everywhere on
$\Omega$. Thus $\nabla S_{\infty} g(u_{\infty})=0$ and in
particular
$$
\dv_y (\nabla S_{\infty} \:u_{\infty}(1-u_{\infty}))=0.
$$

\vskip1mm

\noindent{\it Fourth step. Proof of \eqref{strongcvg}.}

We have already proved that as $k\to\infty$,
$$
\int_0^T\int_\Omega|\nabla S_{n_k}|^2 u_{n_k}
(1-u_{n_k})\to\int_0^T\int_\Omega|\nabla S_{\infty}|^2
u_{\infty}(1-u_{\infty})=0.
$$
The above convergence results entails that
$$
\int_0^T\int_\Omega|\nabla S_{n_k}|^2 u_{n_k}^2
\to\int_0^T\int_\Omega|\nabla S_{\infty}|^2 u_{\infty}^2.
$$
And since $\nabla S_{n_k}\to \nabla S_{\infty}$ in
$L^2((0,T)\times\Omega)$, it follows that
$$
\int_0^T\int_\Omega|\nabla S_{\infty}|^2
\left(u_{n_k}^2-u_{\infty}^2\right) \to 0.
$$
Writing $$(u_{n_k}-u_{\infty})^2=u_{n_k}^2-u_{\infty}^2 - 2
u_{n_k} u_{\infty} + 2 u_{\infty}^2,$$ and using once more the
weak convergence of $u_{n_k}$, we obtain
$$
\int_0^T\int_\Omega|\nabla S|^2 (u_{n_k}-u_{\infty})^2\to 0.
$$
\eqref{strongcvg} follows easily, extracting a further subsequence
if necessary.

\vskip1mm

\noindent{\it Fifth step. Kinetic formulation.}

Let $f=f(t,y,\xi)$ be the weak limit of $\mathbf
1_{\xi<u_{n_k}(t,y)}$. Notice that it is not obvious that $f$ does
not depend on $t$. Then according to the previous steps and to
section \ref{sec:parabolic}, $f$ satisfies
$$
\p_t f + g(\xi)(\xi - S_{\infty}) \p_{\xi} f + g'(\xi) \nabla
S_{\infty} \cdot \nabla f + g'(\xi)(\rho - u_{\infty}f) =\p_{\xi}
m,
$$
where $m$ is a nonnegative measure and $\rho$ is related to $f$ by
equation \eqref{eq:rolimit}. Moreover, since $u_{n_k}(t,y)$
converges to $u_{\infty}(y)$ a.e. on the set $\{y; \nabla
S_{\infty}(y) \neq 0\}$, we deduce that $f(t,y,\xi)=\mathbf
1_{\xi<u_\infty(y)}$ a.e. on $\{\nabla S_{\infty} \neq 0\}$, and thus
$(\rho - u_{\infty}f)=0$ and $g(\xi)\p_{\xi} f=0$ on $\{\nabla
S_{\infty} \neq 0\}$.

\end{proof}

\begin{rem}
In general, stationary states of \eqref{eq:KShyp} are not unique,
even when entropy conditions are required and the mean value on
$\Omega$ is prescribed. Thus, it is not obvious that the whole
sequence $u_k$ should converge to a stationary state $u_{\infty}$.
\end{rem}

%--------------------------------------
\section*{Appendix}

This Appendix is devoted to the rigorous proof of inequality
\eqref{in:comp_f/f2}. Regularizations by convolution are used in
order to justify the nonlinear manipulations which led to equation
\eqref{eq:difference}, as in \cite{artBP,BP}. We focus on the case
when $\Omega$ is an arbitrary bounded domain in $\R^N$, with a
$\mathcal C^1$ boundary, and $\nabla S \cdot n_{\Omega}=0$ on
$\p\Omega$; the case when $\Omega=\Pi_{i=1}^N(0,T_i)$, and $S$
satisfies periodic boundary conditions, is in fact easier, and can
be treated in a similar fashion.

We take $\delta_1,\delta_2>0$ arbitrary, and $\varphi_1\in\mathcal
D(\R)$, $\varphi_2\in\mathcal D(\R^N)$, $\varphi_3\in\mathcal
D(\R)$, with
\begin{gather*}
0\leq \varphi_1,\varphi_1,\varphi_3\leq 1,\\
\int_{\R} \varphi_1=\int_{\R^N} \varphi_2=\int_{\R} \varphi_3=1,\\
\supp \:\varphi_1\subset[-1,0],\qquad
\supp \:\varphi_2\subset B_1,\qquad
\supp \:\varphi_3\subset[-1,1].
\end{gather*}
We set $\delta=(\delta_1,\delta_2)$, and
$$
\varphi_{\delta}(t,y,\xi)=\frac{1}{\delta_1\:\delta_2^{N+1}}\varphi_1\left(\frac
t {\delta_1} \right)\;\varphi_2\left(\frac y {\delta_2}
\right)\:\varphi_3\left(\frac {\xi} {\delta_2} \right),
$$
and for $(t,y,\xi)\in[0,\infty)\times\R^{N+1}$
\begin{gather*}
f_{\delta}(t,y,\xi):=f\ast
\varphi_{\delta}(t,y,\xi)=\int_{\R}\int_{\Omega}\int_{\R}f_{\delta}(t',y',\xi')\varphi_{\delta}(t-t',y-y',\xi-\xi')dt'\:dy'\:d\xi',\\
m_{\delta}:=m\ast \varphi_{\delta}.
\end{gather*}

Then $f_{\delta}$ and $m_{\delta}$ are smooth functions of
$t,y,\xi$ for all $\delta>0$, and $0\leq f_{\delta}\leq 1$,
$m_{\delta}\geq 0$. Moreover, $f_{\delta}$ is a solution of \be
\p_t f_{\delta} + (\xi - S)g(\xi) \p_{\xi} f_{\delta} + g'(\xi)
\nabla_y S \cdot \nabla_y f_{\delta} +
R\ast\varphi_{\delta}=\p_{\xi }m_{\delta} +
r_{\delta},\label{eq:fdelta} \ee and the remainder $r_{\delta}$ is
equal to
\begin{eqnarray*}
r_{\delta}&=&(\xi - S)g(\xi) \p_{\xi} f_{\delta} - \left[(\xi -
S)g(\xi) \p_{\xi} f\right]\ast\varphi_{\delta}\\
&&+ g'(\xi) \nabla_y S \cdot \nabla_y f_{\delta} -\left[ g'(\xi)
\nabla_y S \cdot \nabla_y f\right]\ast\varphi_{\delta}.
\end{eqnarray*}

We wish to stress that equation \eqref{eq:fdelta} holds everywhere
in $(0,\infty)\times \Omega_\delta\times \R$, and not in
$(0,\infty)\times \Omega\times \R$, where
$$
\Omega_\delta:=\{y\in\Omega,\ d(y,\p\Omega)\geq\delta\}.
$$
This yields a small difficulty when integrating equation
\eqref{eq:fdelta} on $(0,\infty)\times \Omega_\delta\times \R$,
because $\nabla S\cdot n_{\Omega_\delta}(y)\neq 0$ on
$\p\Omega_\delta$ even though $\nabla S\cdot n_{\Omega}(y)=0$ on
$\p\Omega$. However, this difficulty can be overcome by using the
regularity of $S$ and of the boundary $\p \Omega$.

Before writing an equation for $f_{\delta}-f_{\delta}^2$, let us
first prove that $r_{\delta}\to 0$ in
$L^1((0,T)\times\Omega\times(-R,R))$ for all $T,R>0$. In the rest
of the appendix, we set $z=(t,y,\xi)\in \R^{N+2}$, with $z_0=t$,
$z_i=y_i$ for $1\leq i\leq N$, $z_{N+1}=\xi$. Accordingly, we
define the differential operators
$$
\p_0=\frac{\p}{\p t},\ \p_i=\frac{\p}{\p y_i}\quad 1\leq i\leq N,\
\p_{N+1}=\frac{\p}{\p \xi}
$$
and we set $Q:=(0,\infty)\times \Omega\times \R$.

Then for instance, we have
\begin{eqnarray*}
&&(\xi - S)g(\xi) \p_{\xi} f_{\delta} - \left[(\xi - S)g(\xi)
\p_{\xi} f\right]\ast\varphi_{\delta}\\
&=&(\xi - S)g(\xi) f\ast\p_{\xi} \varphi_{\delta} - \left[(\xi -
S)g(\xi) f\right]\ast\p_{\xi} \varphi_{\delta}\\
&&+\left[\p_{\xi}\left((\xi - S)g(\xi) \right)
f\right]\ast\varphi_{\delta}\\
&=&\int_{Q}\left(G(z)-G(z')\right)f(z')\p_{N+1}\varphi_{\delta}(z-z')\:dz'\\
&&+\left[\p_{N+1}G\: f\right]\ast\varphi_{\delta}
\end{eqnarray*}
where
$$
G(z):=(\xi - S(t,y))g(\xi),\quad z=(t,y,\xi).
$$
Then $\left[\p_{N+1}G\: f\right]\ast\varphi_{\delta}$ converges to
$\p_{N+1}G\: f$ in $L^p((0,T)\times \Omega \times (R,R))$ for all
$T,R>0$. And setting $\psi_k(z)= z_k \p_{N+1}\varphi(z)$ for
$1\leq k\leq N$, we have
\begin{eqnarray*}
&&\int_{Q}\left(G(z)-G(z')\right)f(z')\p_{N+1}\varphi_{\delta}(z-z')\:dz'\\
&=&\int_{Q}\int_0^1\p_k G(\tau z + (1-\tau)
z')\:f(z')\psi_{k,\delta}\left(z-z' \right)\;dz'\:d\tau
\end{eqnarray*}
The above integral converges to
$$
\p_k G(z) f(z) \int_{\R^{N+2}}\psi_k(z')\:dz'
$$
in $L^2((0,T)\times \Omega \times (-R,R))$ for all $T,R>0$ (recall
that $S$ is bounded in $L^{\infty}(0,T; W^{2,q}(\Omega))\cap
W^{1,2}(0,T; H^1(\Omega))$ thanks to proposition
\ref{prop:epsmodel}). But $\int_{\R^{N+2}}\psi_k(z')\:dz'=0$ if
$k\neq N+1$ and $\int_{\R^{N+2}}\psi_{N+1}(z')\:dz'=-1$. Thus
$$
(\xi - S)g(\xi) \p_{\xi} f_{\delta} - \left[(\xi - S)g(\xi)
\p_{\xi} f\right]\ast\varphi_{\delta}
$$
converges to 0 in $L^2((0,T)\times \Omega \times (-R,R))$ for all
$T,R>0$. The other term can be treated in a similar way, using the
bounds on $S$ derived in proposition \ref{prop:epsmodel}.

We now go back to the equation on $f_{\delta}$; since $f_{\delta}$
is smooth in $t,y,\xi$, we can use the chain rule and write, for
$(t,y,\xi)\in(0,\infty)\times \Omega_\delta\times\R$,
\begin{multline*}
\p_t \left(f_{\delta}-f_{\delta}^2\right) + (\xi - S)g(\xi)
\p_{\xi} \left(f_{\delta}-f_{\delta}^2\right) + g'(\xi) \nabla_y S
\cdot \nabla_y \left(f_{\delta}-f_{\delta}^2\right)+\\+
R\ast\varphi_{\delta}(1-2f_{\delta})=\p_{\xi }m_{\delta}
(1-2f_{\delta})+ r_{\delta}(1-2f_{\delta})
\end{multline*}
We now integrate the above equation on $\Omega_\delta\times\R$;
notice that since $f=0$ for $\xi>1$ and $f=1$ for $\xi<0$, we have
$f_{\delta}-f_{\delta}^2=0$ for $\xi\leq -\delta$ or $\xi\geq 1 +
\delta$, and similarly, $m_{\delta}, r_{\delta}=0$ for $\xi\leq
-\delta$ or $\xi\geq 1 + \delta$. Thus
\begin{eqnarray*}
\frac{d}{dt}\int_{\Omega_\delta\times\R}\left(f_{\delta}-f_{\delta}^2\right)
&=&\int_{\Omega_\delta\times\R}\left(f_{\delta}-f_{\delta}^2\right)\left[\p_{\xi}\left((\xi
- S)g(\xi)\right)+g'(\xi)\Delta S\right]\\
&&-\int_{\Omega_\delta\times\R}R\ast\varphi_{\delta}(1-2f_{\delta})\\
&& +2\int_{\Omega_\delta\times\R}m_{\delta} \p_{\xi}f_{\delta} +
\int_{\Omega_\delta\times\R}r_{\delta}(1-2f_{\delta})\\
&&- \int_{\p\Omega_\delta\times\R}\left(f_{\delta}-f_{\delta}^2\right)g'(\xi)\nabla S \cdot n_{\Omega_\delta}(y)\:dS(y) d\xi\\
&\leq &
C\int_{\Omega_\delta\times\R}\left(f_{\delta}-f_{\delta}^2\right)+C\int_{\Omega_\delta\times\R}\left|f-f^2\right|\ast\varphi_{\delta}\\
&&+ ||r_{\delta}(t)||_{L^1(\Omega\times(-1,2))} + C
\int_{\p\Omega_\delta}\left|\nabla S \cdot n_{\Omega_\delta}(y)
\right|\:dS(y)
\end{eqnarray*}
In the above inequality, we have used the fact that $\p_\xi
f_{\delta}=-\nu\ast\varphi_{\delta}\leq 0$, where $\nu$ was
defined in \eqref{def:nu}, together with lemma \ref{lem:rho_uf}.
Moreover, notice that since the function $x\mapsto x-x^2$ is
concave, by Jensen's inequality, we get
$$
\int_{\Omega_\delta\times\R}\left(f-f^2\right)\ast\varphi_{\delta}\leq
\int_{\Omega_\delta\times\R}\left(f_{\delta}-f_{\delta}^2\right).
$$
And since $S$ belongs to $W^{2,q}$ for all $q<\infty$, $\nabla S
\in \mathcal C^{0,\alpha}(\bar{\Omega})$ for some $0<\alpha<1$;
remember that we have assumed that the boundary $\p\Omega$ is at
least $\mathcal C^1$. In such conditions, it is easily proved that
$$
\int_{\p\Omega_\delta}\left|\nabla S \cdot n_{\Omega_\delta}(y)
\right|\:dS(y)\to 0
$$
as $\delta\to 0$ for almost every $t>0$. In the following, we set
$$
U_{\delta}(t)=||r_{\delta}(t)||_{L^1(\Omega\times(-1,2))} + C
\int_{\p\Omega_\delta}\left|\nabla S \cdot n_{\Omega_\delta}(y)
\right|\:dS(y),
$$
and we have proved that $U_{\delta}\to 0$ in
$L^1_{\text{loc}}([0,\infty))$ as $\delta\to 0$.

Thus we are led to
$$
\frac{d}{dt}\int_{\Omega_\delta\times\R}\left(f_{\delta}-f_{\delta}^2\right)\leq
C\int_{\Omega_\delta\times\R}\left(f_{\delta}-f_{\delta}^2\right)
+U_{\delta}(t).
$$

Consequently, by Gronwall's lemma,
\be\int_{\Omega_\delta\times\R}\left(f_{\delta}(t)-f_{\delta}(t)^2\right)\leq
e^{Ct}\int_{\Omega_\delta\times\R}\left(f_{\delta}(t=0)-f_{\delta}(t=0)^2\right)
+ \int_0^t e^{C(t-s)}U_{\delta}(s)\:ds.\label{in:Gronwall}\ee

There only remains to prove that
$f_{\delta}(t=0)-f_{\delta}(t=0)^2$ goes to 0 as $\delta\to 0$.
This is a consequence of the fact that $f_\delta(t=0)$ strongly
converges to $f(t=0)=\mathbf 1_{\xi<u_0}$, but the latter is not
obvious since
$$
f_\delta(t=0,y,\xi)=\int_0^{\infty}\int_{\Omega\times\R}f(t',y',\xi')\varphi_{\delta}(-t',y-y',\xi-\xi')\;dt'\:dy'\:d\xi'.
$$

We therefore use the same technique as in \cite{BP} : since the
proof is strictly identical to the one in \cite{BP}, we only
recall briefly the main arguments.

Let
$$
T_{\delta_1}(t)=1-\frac{1}{\delta_1}\int_0^t\varphi_1\left(-\frac
s {\delta_1}\right)\:ds;
$$
then $\p_t
T_{\delta_1}(t)=-\frac{1}{\delta_1}\varphi_1\left(-\frac
t{\delta_1}\right)$, and thus for
$(t,y,\xi)\in(0,\infty)\times\Omega_\delta\times\R$,
$f_\delta(t=0,y,\xi)$ can be written as
\begin{eqnarray*}
f_\delta(t=0,y,\xi)&=&I_{\delta}(y,\xi)\\
&&+
\p_\xi\int_0^\infty\!\!\!\int_{\R^{N+1}}m(t',y',\xi')\frac{1}{\delta_2^{N+1}}\varphi_2\left(\frac{y-y'}{\delta_2}
\right)\varphi_3\left(\frac{\xi-\xi'}{\delta_2}
\right)\;dt'\:dy'\:d\xi'\\
&&+\int_{\R^{N+1}}\mathbf
1_{\xi<u_0(y)}\frac{1}{\delta_2^{N+1}}\varphi_2\left(\frac{y-y'}{\delta_2}
\right)\varphi_3\left(\frac{\xi-\xi'}{\delta_2}
\right)\:dy'\:d\xi'
\end{eqnarray*}
where $||I_{\delta}||_{L^{\infty}}\leq C \delta_1/\delta_2$.
Passing first to the weak limit as $\delta_1,\delta_2 \to 0$ with
$\delta_1/\delta_2\to 0$ in the above equation entails that the
weak limit $F$ of $f_\delta(t=0,y,\xi)$ satisfies
$$
F=\p_{\xi} M + 1_{\xi<u_0(y)}
$$
for some nonnegative measure $M$ vanishing for large $\xi$. This
leads to $F=1_{\xi<u_0(y)}$ and $M=0$ thanks to a lemma in
\cite{BP}. Then, the above formula for $f_\delta(t=0,y,\xi)$ is
used once again to find the weak limit of $f_\delta(t=0,y,\xi)^2$.
It is easily proved that
$$
f_\delta(t=0,y,\xi)^2\rightharpoonup F^2=1_{\xi<u_0(y)},
$$
and thus the convergence is strong.

Consequently, $f_\delta-f_\delta^2$ converges to 0 in
$L^1_\text{loc}((0,\infty)\times\Omega\times\R)$. Since
$f_\delta-f_\delta^2\to f-f^2$ in $L^1_\text{loc}(0,\infty;
L^1(\Omega\times\R))$, we deduce that $f=0$ or $f=1$ almost
everywhere. The rest of the proof, exposed in section
\ref{sec:rigidity}, is therefore justified.

\bibliography{articlesKS}
\end{document}